\documentclass{amsart}%
\usepackage{amssymb}
\usepackage{amsfonts}
\usepackage{amsmath}
\usepackage{graphicx}%
\newtheorem{theorem}{Theorem}
\theoremstyle{plain}

\newtheorem{lemma}{Lemma}

\numberwithin{equation}{section}
\begin{document}
\title[ ]{On a Sharp Inequality Relating Yamabe Invariants on a Poincare-Einstein Manifold}
\author{Xiaodong Wang}
\address{Department of Mathematics, Michigan State University, East Lansing, MI 48824}
\email{xwang@msu.edu}
\author{Zhixin Wang}
\address{Department of Mathematics, Michigan State University, East Lansing, MI 48824}
\email{wangz117@msu.edu}

\begin{abstract}
For a Poincare-Einstein manifold under certain restrictions, X. Chen, M. Lai
and F. Wang \cite{CLW} proved a sharp inequality relating Yamabe invariants.
We show that the inequality is true without any restriction.

\end{abstract}
\maketitle

A complete Riemannian manifold $\left(  X^{n},g_{+}\right)  $ is called
conformally compact if $X$ is the interior of a compact manifold $\overline
{X}$ with nonempty boundary $\Sigma=\partial\overline{X}$ and for any defining
function $r$ on $\overline{X}$ the metric $\overline{g}:=r^{2}g_{+}$ extends
to a smooth metric on $\overline{X}$ (we gloss over the boundary regularity of
$\overline{g}$ as it does not play an essential role in our discussion). In
particular, $\overline{\gamma}=\overline{g}|_{\partial\overline{X}}$ defines a
Riemannian metric on $\Sigma$. Though a different choice of $r$ gives rise to
a different $\overline{g}$, the conformal classes $\left[  \overline
{g}\right]  $ and $\left[  \overline{g}|_{\Sigma}\right]  $ are invariantly
defined. $\left(  \Sigma,\left[  \overline{g}|_{\Sigma}\right]  \right)  $ is
called the conformal infinity of $\left(  X^{n},g_{+}\right)  $. If
furthermore $Ric\left(  g_{+}\right)  =-\left(  n-1\right)  g_{+}$, then
$\left(  X^{n},g_{+}\right)  $ is called Poincar\'{e}-Einstein (or conformally
compact Einstein). Poincar\'{e}-Einstein manifolds have been studied
intensively since the work of Graham-Lee \cite{GL}. More recently the subject
has gained new momentum from string theory or more specifically the AdS/CFT
correspondence in which such manifolds serve as the framework for a connection
between supergravity and conformal field theory. 

A guiding principle for studying Poincar\'{e}-Einstein manifolds is to
understand the interaction between the geometry on $X$ and the conformal
geometry on $\partial\overline{X}$. Recently Gursky-Han proved an eminent
result of this type which relates the Yamabe invariant of $\left(
\overline{X},\left[  \overline{g}\right]  \right)  $ and the Yamabe invariant
of $\left(  \Sigma,\left[  \overline{g}|_{\Sigma}\right]  \right)  $. Before
stating their result, let us first recall some basics on the Yamabe problem.

On a compact Riemannian manifold $\left(  M^{n},g\right)  $ with nonempty
boundary $\Sigma=\partial M$, there are in fact two types of Yamabe problems.
We first consider the functional%
\[
E_{g}\left(  u\right)  =\int_{M}\left(  \frac{4\left(  n-1\right)  }%
{n-2}\left\vert \nabla u\right\vert ^{2}+Ru^{2}\right)  dv_{g}+2\int_{\Sigma
}Hu^{2}d\sigma_{g},
\]
where $R$ is the scalar curvature and $H$ is the mean curvature of the
boundary. (In our convention $H$ is the trace of the 2nd fundamental form.)
This functional has the important property of being conformally invariant: if
$\widetilde{g}=\phi^{4/\left(  n-2\right)  }g$ is another metric, then
$E_{\widetilde{g}}\left(  u\right)  =E_{g}\left(  u\phi\right)  $. The type I
Yamabe invariant is defined as
\[
Y\left(  M,\left[  g\right]  \right)  =\inf_{u\in H^{1}\left(  M\right)
\backslash\left\{  0\right\}  }\frac{E_{g}\left(  u\right)  }{\left(  \int
_{M}\left\vert u\right\vert ^{2n/\left(  n-2\right)  }dv_{g}\right)  ^{\left(
n-2\right)  /n}}.
\]
As the notation indicates, $Y\left(  M,\left[  g\right]  \right)  $ depends
only on the conformal class $\left[  g\right]  $. If the infimum is achieved,
then there exists an essentially unique minimizer $u$ which is smooth and
positive and the metric $u^{4/\left(  n-2\right)  }g$ then has constant scalar
curvature on $M$ and zero mean curvature on $\Sigma$. The type I Yamabe
problem whether the infimum is always achieved has not been completely solved.
It has been solved in many cases (see \cite{E1} and \cite{BC}).

The type II Yamabe invariant is defined as
\[
Q\left(  M,\Sigma,\left[  g\right]  \right)  =\inf_{u\in H^{1}\left(
M\right)  \backslash\left\{  0\right\}  }\frac{E_{g}\left(  u\right)
}{\left(  \int_{\Sigma}\left\vert u\right\vert ^{2\left(  n-1\right)  /\left(
n-2\right)  }d\sigma_{g}\right)  ^{\left(  n-2\right)  /\left(  n-1\right)  }%
}.
\]
It should be noted that $Q\left(  M,\Sigma,\left[  g\right]  \right)  $ can be
$-\infty$. If $Q\left(  M,\Sigma,\left[  g\right]  \right)  $ $>-\infty$ and
the infimum is achieved, then there exists an essentially unique minimizer $u$
which is smooth and positive and the metric $u^{4/\left(  n-2\right)  }g$ then
has zero scalar curvature on $M$ and constant mean curvature on $\Sigma$. The
type II Yamabe problem whether the infimum is achieved when $Q\left(
M,\Sigma,\left[  g\right]  \right)  $ $>-\infty$ has been solved in various
cases (see \cite{E2},\cite{M1} and \cite{M2}). But there are still cases that
remain open. Apart from the minimization problem, both $Y\left(  M,\left[
g\right]  \right)  $ and $Q\left(  M,\Sigma,\left[  g\right]  \right)  $ are
important invariants for which it is useful to have good estimates.

We now come back to Poincar\'{e}--Einstein manifolds. Let $\left(  X^{n}%
,g_{+}\right)  $ be a Poincar\'{e}--Einstein manifold and $\Sigma
=\partial\overline{X}$. We pick a fixed defining function $r$ on $\overline
{X}$ which gives rise to a metric $\overline{g}=r^{2}g_{+}$ on $\overline{X}$.
As $\left[  \overline{g}\right]  $ and $\left[  \overline{g}|_{\Sigma}\right]
$ are invariantly defined, the Yamabe invariants $Y\left(  \overline
{X},\left[  \overline{g}\right]  \right)  ,Q\left(  \overline{X}%
,\Sigma,\left[  \overline{g}\right]  \right)  $ and $Y\left(  \Sigma,\left[
\overline{g}|_{\Sigma}\right]  \right)  $ are natural invariants of $\left(
X^{n},g_{+}\right)  $. We can now state the result of Gursky-Han.

\begin{theorem}
(Gursky-Han \cite{GH}) Suppose $\overline{X}$ satisfies one of the following
two conditions

\begin{enumerate}
\item $3\leq n\leq5$;

\item $n\geq6$ and $M$ is spin.
\end{enumerate}

Let $\widetilde{g}$ be a type I Yamabe minimizer in $\left[  \overline
{g}\right]  $. Then%
\begin{align*}
\frac{n}{n-2}Y\left(  \Sigma,\left[  \overline{g}|_{\Sigma}\right]  \right)
&  \leq\frac{n-2}{4\left(  n-1\right)  }Y\left(  \overline{X},\left[
\overline{g}\right]  \right)  I^{2},\text{ if }n\geq4;\\
12\pi\chi\left(  \Sigma\right)   &  \leq\frac{n-2}{4\left(  n-1\right)
}Y\left(  \overline{X},\left[  \overline{g}\right]  \right)  I^{2},\text{ if
}n=3,
\end{align*}
where $I$ is the isoperimetric ratio for $\widetilde{g}$ (cf. \cite{GH} for
the precise definition). Moreover, if the equality holds, then $\widetilde{g}$
is Einstein and $\widetilde{g}|_{\Sigma}$ has constant scalar curvature.
\end{theorem}

Modifying Gursky-Han's method, X. Chen, M. Lai and F. Wang \cite{CLW} proved
the following result for type II Yamabe invariant.

\begin{theorem}
(Chen-Lai-Wang \cite{CLW}) Let $\left(  X^{n},g_{+}\right)  $ be a
Poincar\'{e}--Einstein manifold s.t. $\left(  \overline{X},\overline
{g}\right)  $ satisfies one of the conditions

\begin{enumerate}
\item[(a)] the dimension $3\leq n\leq7$;

\item[(b)] the dimension $n\geq8$ and $M$ is spin;

\item[(c)] the dimension $n\geq8$ and $M$ is locally conformally flat.
\end{enumerate}

Then
\begin{align*}
Y\left(  \Sigma,\left[  \overline{g}|_{\Sigma}\right]  \right)   &  \leq
\frac{n-2}{4\left(  n-1\right)  }Q\left(  \overline{X},\Sigma,\left[
\overline{g}\right]  \right)  ^{2},\text{ if }n\geq4;\\
32\pi\chi\left(  \partial\overline{X}\right)   &  \leq Q\left(  \overline
{X},\Sigma,\left[  \overline{g}\right]  \right)  ^{2},\text{ if }n=3.
\end{align*}
Moreover, the equality holds if and only if $\left(  X^{n},g_{+}\right)  $ is
isometric to the hyperbolic space $\left(  \mathbb{H}^{n},g_{\mathbb{H}%
}\right)  $.
\end{theorem}

The inequality in Theorem 2 is very elegant as both sides are natural
invariants. But the proof makes use of a type II Yamabe minimizer which exist
under the extra conditions. Due to the fact that the type II Yamabe problem is
not completely solved yet, the result is not proved for all
Poincar\'{e}-Einstein manifolds.

It should be noted that Theorem 2 without any restriction is claimed in
\cite{R}. In the exceptional cases not covered by Theorem 2, he applies the
same method in \cite{CLW} to a metric in $\left[  \overline{g}\right]  $ with
zero scalar curvature on $\overline{X}$ and constant mean curvature on
$\Sigma$, whose existence is proved in \cite{MN}. But it is clear that this
metric is in general just a critical point of the Yamabe functional not a
minimizer. As such the argument does not yield any information on the Yamabe
invariant. Therefore the proof in \cite{R} for the exceptional cases is invalid.

In this note, we remove the restrictions in Theorem 2 by an indirect route.
Note the inequality is vacuous when $Y\left(  \partial X,\left[  g\right]
\right)  \leq0$. Therefore we state the result in the following way.

\begin{theorem}
Let $\left(  X^{n},g_{+}\right)  $ be a Poincar\'{e}--Einstein manifold whose
conformal infinity has nonnegative Yamabe invariant. Then%
\begin{align*}
Q\left(  \overline{X},\Sigma,\left[  \overline{g}\right]  \right)   &
\geq2\sqrt{\frac{\left(  n-1\right)  }{\left(  n-2\right)  }Y\left(
\Sigma,\left[  \overline{g}|_{\Sigma}\right]  \right)  }\text{ if }n\geq4;\\
Q\left(  \overline{X},\Sigma,\left[  \overline{g}\right]  \right)   &
\geq4\sqrt{2\pi\chi\left(  \Sigma\right)  }\text{ if }n=3.
\end{align*}
Moreover, the equality holds iff $\left(  X^{n},g_{+}\right)  $ is isometric
to the hyperbolic space $\left(  \mathbb{H}^{n},g_{\mathbb{H}}\right)  $.
\end{theorem}

We remark that by \cite{WY} and \cite{CG} $\Sigma$ is connected when it has
nonnegative Yamabe invariant. We first address the equality case.  First
suppose $n\geq4$. If $Y\left(  \Sigma,\left[  \overline{g}|_{\Sigma}\right]
\right)  <Y\left(  \mathbb{S}^{n-1},g_{\mathbb{S}^{n-1}}\right)  $, then the
equality implies
\[
Q\left(  \overline{X},\partial\overline{X},\left[  \overline{g}\right]
\right)  <Q\left(  \overline{\mathbb{B}^{n}},\partial\overline{\mathbb{B}^{n}%
},\left[  g_{0}\right]  \right)  ,
\]
the Yamabe invariant of the model case. It is well known that $Q\left(
\overline{X},\Sigma,\left[  \overline{g}\right]  \right)  $ is then achieved
and therefore the rigidity is covered by the original argument in \cite{CLW}.
If $Y\left(  \Sigma,\left[  \overline{g}|_{\Sigma}\right]  \right)  =Y\left(
\mathbb{S}^{n-1},g_{\mathbb{S}^{n-1}}\right)  $, then by the solution of the
Yamabe problem on closed manifolds, $\left(  \Sigma,\left[  \overline
{g}|_{\Sigma}\right]  \right)  $ is conformally diffeomorphic to
$\mathbb{S}^{n-1}$. Then $\left(  X^{n},g_{+}\right)  $ is isometric to the
hyperbolic space $\left(  \mathbb{H}^{n},g_{\mathbb{H}}\right)  $ by \cite{DJ}
and \cite{LQS}. The same argument works when $n=3$.

To prove the inequality, we first recall some basic facts on
Poincar\'{e}-Einstein manifolds. Let $h\in\left[  \overline{g}|_{\Sigma
}\right]  $ be a metric on $\Sigma$. It is proved in \cite{Lee} that there is
a defining function $r$ s.t. in a collar neighborhood of $\Sigma$%
\[
g_{+}=r^{-2}\left(  dr^{2}+h_{r}\right)  ,
\]
where $h_{r}$ is an $r$-dependent family of metrics on $\partial\overline{X}$
with $h_{r}|_{r=0}=h$. Moreover we have the following expansion (see, e.g.
\cite{GW})%
\[
h_{r}=h+h_{2}r^{2}+o\left(  r^{2}\right)  ,
\]
where%
\[
h_{2}=\left\{
\begin{array}
[c]{cc}%
-\frac{1}{n-3}\left(  Ric\left(  h\right)  -\frac{R_{h}}{2\left(  n-2\right)
}h\right)  , & \text{if }n\geq4;\\
-\frac{1}{4}h, & \text{if }n=3.
\end{array}
\right.
\]
It follows that $\overline{g}=r^{2}g_{+}$ has totally geodesic boundary. As we
assume $Y\left(  \Sigma,\left[  \overline{g}|_{\Sigma}\right]  \right)  \geq
0$, we can choose $h$ to have $R_{h}\geq0$.

Lee \cite{Lee} constructed a positive smooth function $\phi$ on $X$ s.t.
$\Delta\phi=n\phi$ and near $\partial\overline{X}$%
\[
\phi=r^{-1}+\frac{R_{h}}{4\left(  n-1\right)  \left(  n-2\right)  }r+O\left(
r^{2}\right)  .
\]
Under the condition $R_{h}\geq0$, he further proved that $\left\vert
d\phi\right\vert _{g_{+}}^{2}\leq\phi^{2}$. It was observed by Qing \cite{Q}
that $\widetilde{g}:=\phi^{-2}g_{+}$ is then a metric on $\overline{X}$ with
nonnegative scalar curvature and totally geodesic boundary and hence
$E_{\widetilde{g}}\left(  u\right)  \geq0$. By the conformal invariance, we
also have $E_{\overline{g}}\left(  u\right)  \geq0$. For $1<1\leq n/\left(
n-2\right)  $, consider%
\[
\lambda_{q}:=\inf\frac{E_{\overline{g}}\left(  u\right)  }{\left(
\int_{\Sigma}\left\vert u\right\vert ^{q+1}d\sigma_{\overline{g}}\right)
^{2/\left(  q+1\right)  }}.
\]
Notice that $\lambda_{q}=Q\left(  X,\partial X,\left[  \overline{g}\right]
\right)  $ when $q=n/\left(  n-2\right)  $. As $E_{\overline{g}}\left(
u\right)  \geq0$, it is easy to see that $\lim_{q\rightarrow n/\left(
n-2\right)  }\lambda_{q}=Q\left(  X,\partial X,\left[  \overline{g}\right]
\right)  $. Therefore Theorem follows from the following

\begin{theorem}
\label{isub}Let $\left(  X^{n},g_{+}\right)  $ be a Poincar\'{e}--Einstein
manifold whose conformal infinity has positive Yamabe invariant. For
$1<q<n/\left(  n-2\right)  $ the invariant $\lambda_{q}$ satisfies%
\begin{align*}
\lambda_{q}  &  \geq2\sqrt{\frac{\left(  n-1\right)  }{\left(  n-2\right)
}Y\left(  \Sigma,\left[  \overline{g}|_{\Sigma}\right]  \right)  }V\left(
\Sigma,\overline{g}\right)  ^{-\frac{\left(  n-q\left(  n-2\right)  \right)
}{\left(  n-3\right)  \left(  q+1\right)  }}\text{ if }n\geq4;\\
\lambda_{q}  &  \geq4\sqrt{2\pi\chi\left(  \Sigma\right)  }V\left(
\Sigma,\overline{g}\right)  ^{-\frac{3-q}{2\left(  q+1\right)  }}\text{ if
}n=3.
\end{align*}

\end{theorem}

The rest of the note is devoted to the proof of this theorem. We recall the
following lemma from \cite{CLW} which follows by a direct calculation.

\begin{lemma}
Suppose $g\in\left[  \overline{g}\right]  $ is scalar flat. Near
$\Sigma=\partial\overline{X}$, we can write
\[
g=dr^{2}+g_{ij}\left(  r,x\right)  dx_{i}dx_{j},
\]
where $\left\{  x_{1},\cdots,x_{n-1}\right\}  $ are local coordinates on
$\Sigma$. If we write $g_{+}=\rho^{-2}g$, then
\[
\rho=r-\frac{H}{2\left(  n-1\right)  }r^{2}+\frac{1}{6}\left(  \frac
{R^{\Sigma}}{n-2}-\frac{H^{2}}{n-1}\right)  r^{3}+o\left(  r^{3}\right)  .
\]
In particular,%
\[
\rho^{-1}\left(  \frac{\partial}{\partial\nu}\left\vert \nabla\rho\right\vert
^{2}+\rho^{-1}\left(  1-\left\vert \nabla\rho\right\vert ^{2}\right)
\frac{\partial\rho}{\partial\nu}\right)  |_{\Sigma}=\frac{R^{\Sigma}}%
{n-2}-\frac{H^{2}}{n-1}.
\]

\end{lemma}

Since the trace operator $H^{1}\left(  \overline{X}\right)  \rightarrow
L^{q+1}\left(  \Sigma\right)  $ is compact for $q<n/\left(  n-2\right)  $, by
standard elliptic theory, the above infimum $\lambda_{q}$ is achieved by a
smooth, positive function $u$ s.t.%
\begin{equation}
\int_{\Sigma}u^{q+1}d\overline{\sigma}=1\label{Norm}%
\end{equation}
and%
\begin{align*}
-\frac{4\left(  n-1\right)  }{n-2}\overline{\Delta}u+\overline{R}u &  =0\text{
on }\overline{X},\\
\frac{4\left(  n-1\right)  }{n-2}\frac{\partial u}{\partial\overline{\nu}%
}+2\overline{H}u &  =\lambda_{q}u^{q}\text{ on }\Sigma.
\end{align*}
The equations imply that the metric $g=u^{4/\left(  n-2\right)  }\overline{g}$
has zero scalar curvature and the mean curvature on the boundary is given by
\begin{equation}
H=\frac{\lambda_{q}}{2}u^{q-\frac{n}{n-2}}.\label{Hf}%
\end{equation}

We now apply the method of \cite{GH} and \cite{CLW} to $g$. Write $g_{+}%
=\rho^{-2}g$. As $g_{+}$ is Einstein, we have
\[
E_{g}=-\left(  n-2\right)  \rho^{-1}\left[  D^{2}\rho-\frac{\Delta\rho}%
{n}g\right]  .
\]
Let $r$ be the distance function to $\Sigma$ w.r.t. $g$ and denote
$X_{\varepsilon}=\left\{  r\geq\varepsilon\right\}  $. From the above
identity, we have%
\begin{align*}
\int_{X_{\varepsilon}}\rho\left\vert E_{g}\right\vert ^{2}dv_{g} &  =-\left(
n-2\right)  \int_{X_{\varepsilon}}\left\langle D^{2}\rho-\frac{\Delta\rho}%
{n}g,E_{g}\right\rangle dv_{g}\\
&  =-\left(  n-2\right)  \int_{X_{\varepsilon}}\left\langle D^{2}\rho
,E_{g}\right\rangle dv_{g}\\
&  =-\left(  n-2\right)  \int_{\partial X_{\varepsilon}}E\left(  \nabla
\rho,\nu\right)  dv_{g}\\
&  =\left(  n-2\right)  ^{2}\int_{\partial X_{\varepsilon}}\rho^{-1}\left(
D^{2}\rho\left(  \nabla\rho,\nu\right)  -\frac{\Delta\rho}{n}\frac
{\partial\rho}{\partial\nu}\right)  dv_{g}\\
&  =\frac{\left(  n-2\right)  ^{2}}{2}\int_{\partial X_{\varepsilon}}\rho
^{-1}\left(  \frac{\partial}{\partial\nu}\left\vert \nabla\rho\right\vert
^{2}-\frac{2\Delta\rho}{n}\frac{\partial\rho}{\partial\nu}\right)  dv_{g}.
\end{align*}
As $g$ has zero scalar curvature and $g_{+}=\rho^{-2}g$ has scalar curvature
$-n\left(  n-1\right)  $, we have%
\[
-\frac{2}{n}\rho\Delta\rho=1-\left\vert \nabla\rho\right\vert ^{2}.
\]
Thus
\[
\int_{X_{\varepsilon}}\rho\left\vert E_{g}\right\vert ^{2}dv_{g}=\frac{\left(
n-2\right)  ^{2}}{2}\int_{\partial X_{\varepsilon}}\rho^{-1}\left(
\frac{\partial}{\partial\nu}\left\vert \nabla\rho\right\vert ^{2}+\rho
^{-1}\left(  1-\left\vert \nabla\rho\right\vert ^{2}\right)  \frac
{\partial\rho}{\partial\nu}\right)  dv_{g}.
\]
Note that $\nu=-\frac{\partial}{\partial r}$ along $\partial X_{\varepsilon}$.
Letting $\varepsilon\rightarrow0$, by Lemma 1 we obtain%
\begin{equation}
\frac{2}{\left(  n-2\right)  ^{2}}\int_{X}\rho\left\vert E_{g}\right\vert
^{2}dv_{g}=\int_{\Sigma}\left(  \frac{H^{2}}{n-1}-\frac{R^{\Sigma}}%
{n-2}\right)  d\sigma.\label{fe1}%
\end{equation}

\bigskip By (\ref{Norm}) and the Holder inequality again%
\begin{align*}
\int_{\Sigma}H^{2}d\sigma &  =\left(  \frac{\lambda_{q}}{2}\right)  ^{2}%
\int_{\Sigma}u^{2\left(  q-\frac{n}{n-2}\right)  }u^{2\left(  n-1\right)
/\left(  n-2\right)  }d\overline{\sigma}\\
&  =\left(  \frac{\lambda_{q}}{2}\right)  ^{2}\int_{\Sigma}u^{2\left(
q-\frac{1}{n-2}\right)  }d\overline{\sigma}\\
&  \leq\left(  \frac{\lambda_{q}}{2}\right)  ^{2}\left(  \int_{\Sigma}%
u^{q+1}d\overline{\sigma}\right)  ^{2\left(  q-\frac{1}{n-2}\right)  /\left(
q+1\right)  }V\left(  \Sigma,\overline{g}\right)  ^{\left(  \frac{n}%
{n-2}-q\right)  /\left(  q+1\right)  }\\
&  =\left(  \frac{\lambda_{q}}{2}\right)  ^{2}V\left(  \Sigma,\overline
{g}\right)  ^{\left(  \frac{n}{n-2}-q\right)  /\left(  q+1\right)  }.
\end{align*}
Plugging the above inequality into (\ref{fe1}), we obtain%
\begin{equation}
\frac{2}{\left(  n-2\right)  ^{2}}\int_{X}\rho\left\vert E_{g}\right\vert
^{2}dv_{g}\leq\frac{\lambda_{q}^{2}}{4\left(  n-1\right)  }V\left(
\Sigma,\overline{g}\right)  ^{\left(  \frac{n}{n-2}-q\right)  /\left(
q+1\right)  }-\frac{1}{n-2}\int_{\Sigma}R^{\Sigma}d\sigma.\label{fe2}%
\end{equation}
When $n=3$, this implies%
\[
\lambda_{q}^{2}V\left(  \Sigma,\overline{g}\right)  ^{\left(  3-q\right)
/\left(  q+1\right)  }\geq32\pi\chi\left(  \Sigma\right)  .
\]
In the following, we assume $n>3\,$. By (\ref{Norm}) and the Holder
inequality
\begin{align*}
1 &  =\int_{\Sigma}u^{q+1}d\overline{\sigma}\\
&  \leq\left(  \int_{\Sigma}u^{2\left(  n-1\right)  /\left(  n-2\right)
}d\overline{\sigma}\right)  ^{\frac{\left(  q+1\right)  \left(  n-2\right)
}{2\left(  n-1\right)  }}V\left(  \Sigma,\overline{g}\right)  ^{\frac
{n-q\left(  n-2\right)  }{2\left(  n-1\right)  }}\\
&  =V\left(  \Sigma,g\right)  ^{\frac{\left(  q+1\right)  \left(  n-2\right)
}{2\left(  n-1\right)  }}V\left(  \Sigma,\overline{g}\right)  ^{\frac
{n-q\left(  n-2\right)  }{2\left(  n-1\right)  }}%
\end{align*}
Thus%
\[
V\left(  \Sigma,\overline{g}\right)  ^{-\frac{n-q\left(  n-2\right)  }{\left(
n-2\right)  \left(  q+1\right)  }}\leq V\left(  \Sigma,g\right)  .
\]
Plugging this inequality into (\ref{fe2}) yields
\begin{align*}
\frac{2}{\left(  n-2\right)  ^{2}}\int_{M}\rho\left\vert E_{g}\right\vert
^{2}dv_{g} &  \leq\frac{V\left(  \Sigma,g\right)  ^{\frac{n-1}{n-3}}}{4\left(
n-1\right)  }\left[  \lambda_{q}^{2}V\left(  \Sigma,\overline{g}\right)
^{\frac{2\left(  n-q\left(  n-2\right)  \right)  }{\left(  n-3\right)  \left(
q+1\right)  }}-\frac{4\left(  n-1\right)  }{\left(  n-2\right)  V\left(
\Sigma,g\right)  ^{\frac{n-1}{n-3}}}\int_{\Sigma}R^{\Sigma}d\sigma\right]  \\
&  \leq\frac{V\left(  \Sigma,g\right)  ^{\frac{n-1}{n-3}}}{4\left(
n-1\right)  }\left[  \lambda_{q}^{2}V\left(  \Sigma,\overline{g}\right)
^{\frac{2\left(  n-q\left(  n-2\right)  \right)  }{\left(  n-3\right)  \left(
q+1\right)  }}-\frac{4\left(  n-1\right)  }{\left(  n-2\right)  }Y\left(
\Sigma,\left[  \gamma\right]  \right)  \right]  .
\end{align*}
Therefore
\[
\lambda_{q}^{2}\geq\frac{4\left(  n-1\right)  }{\left(  n-2\right)  }Y\left(
\Sigma\right)  V\left(  \Sigma,\overline{g}\right)  ^{-\frac{2\left(
n-q\left(  n-2\right)  \right)  }{\left(  n-3\right)  \left(  q+1\right)  }}.
\]
This finishes the proof of Theorem \ref{isub}

\end{document}